\documentclass[11pt,onecolumn]{article}													% paper
\usepackage{etoolbox,booktabs}
\usepackage{algorithm, algorithmic, setspace}
\usepackage{graphicx} % for pdf, bitmapped graphics files
\usepackage{amsmath} % assumes amsmath package installed
\usepackage{amssymb}  % assumes amsmath package installed
\usepackage{amsthm}
\usepackage{amsfonts}
\usepackage{color}
\usepackage[english]{babel}

\newcommand{\argmin}[1]{\underset{#1}{\mathrm{argmin\,}}}

\newcommand{\xtrue}{\widetilde{x}}

\newcommand{\sign}{\text{sign}}
\newcommand{\xstar}{x^{\star}}
\newcommand{\lstar}{\lambda^{\star}}
\newcommand{\R}{\mathbb{R}}

\newcommand{\I}{\mathrm{1}}

\newcommand{\soft}{\mathrm{S}}

\newcommand{\epsi}{\epsilon}

\newtheorem{lemma}{Lemma}
\newtheorem{proposition}{Proposition}

\newtheorem{remark}{Remark}

\newtheorem{assumption}{Assumption}
\newtheorem{problem}{Problem}

\title{Integral control of the proximal gradient method\\ for unbiased sparse optimization}
\author{V. Cerone, S. M. Fosson, A. Re, D. Regruto
\thanks{The authors are with the Department of Control and Computer Engineering, Politecnico di Torino, Italy;
    e-mail: sophie.fosson@polito.it.  Funded by the European Union - NextGenerationEU, Mission 4 Component
1.5 - ECS00000036 - CUP E13B22000020001.}}
\begin{document}
\maketitle
\thispagestyle{empty}
\pagestyle{empty}
\begin{abstract} 
Proximal gradient methods are popular in sparse optimization as they are straightforward to implement. Nevertheless, they achieve biased solutions, requiring many iterations to converge. This work addresses these issues through a suitable feedback control of the algorithm's hyperparameter. Specifically, by designing an integral control that does not substantially impact the computational complexity, we can reach an unbiased solution in a reasonable number of iterations. In the paper, we develop and analyze the convergence of the proposed approach for strongly-convex problems. Moreover, numerical simulations validate and extend the theoretical results to the non-strongly convex framework.
\end{abstract}

%\begin{IEEEkeywords} 
%Wireless networked control systems, nonlinear systems, power control
%\end{IEEEkeywords}

%=====================
%
%
\section{Introduction}\label{sec:IN}
Nowadays, parsimonious models, i.e., models that depend on a relatively small number of parameters, play a central role in machine learning, system identification and neural networks. While over-parametrization may have benefits in deep networks \cite{pil25}, parsimony is crucial for diverse purposes: it reduces the computational complexity for lightweight implementation, e.g., in mobile and cyber-physical applications; it provides interpretable representations of physical dynamical systems by selecting the most relevant variables; it prevents overfitting; it deals with compressed measurements and missing data. We refer the reader to, e.g., \cite{bru19book,has15book,fou13} for a comprehensive overview.

In most cases, building parsimonious models from data consists in finding sparse solutions (i.e., solutions with many zeros) to minimization problems, which we refer to as {\it{sparse optimization}}. A valuable approach to promote sparsity and select the most important features is to add a regularization to the cost function to minimize. In the literature, considerable attention is devoted to $\ell_1$ regularization, since the $\ell_1$ norm is the best convex approximation of the number of non-zero components of a vector; see, e.g., \cite{tib96}.

Since sparsity-promoting regularization is usually non-differentiable, the proximal gradient method (PGM) is the natural alternative to gradient descent methods. PGM is an iterative algorithm that consists of a gradient step over the (differentiable) cost function and a proximal operator over the regularization; see, e.g., \cite{com11,par13,fox23} for details. In case of $\ell_1$ regularization, PGM is also known as iterative shrinkage-thresholding algorithm (ISTA, \cite{dau04}), as the proximal map of the $\ell_1$ norm shrinks and thresholds the current estimate.

A drawback of sparsity-promoting regularization is its inherently biased solution. The regularized problem calls for a tradeoff between minimization of the cost function and sparsity: usually, a satisfactory variable selection comes with an unavoidable inaccuracy in assessing the values of the selected variables.

The literature has devoted much attention to this issue for the Lasso estimator, i.e., for the $\ell_1$-regularized least-squares minimization \cite{tib96}. In particular, several works focus on non-convex regularization to correct the bias of Lasso; see, e.g., \cite{zha10MCP,woo16,rak19,fox21sysid}.
Moreover, in \cite{fox23}, the authors notice that PGM applied to Lasso with non-convex regularization yields a faster convergence with respect to ISTA. %The rationale behind this acceleration is as follows. The weight $\lambda>0$ of the $\ell_1$ regularization of Lasso turns out to be the shrinkage-tresholding hyperparameter of ISTA. A large $\lambda$ yields more sparsity, and consequently more bias, while it also speeds up ISTA.
This approach, denoted as AD-ISTA, has an adaptive shrinkage hyperparameter that speeds up the convergence while mitigating the bias effect. From a feedback control perspective, AD-ISTA is a discrete-time dynamical system whose shrinkage hyperparameter is a control input that evolves as a function of the current state value.

Starting from this general feedback perspective, this work proposes a novel approach to tune the hyperparameter of ISTA to keep the velocity of the adaptive approach while optimizing the bias control. More precisely, we investigate a control-theoretic approach and we design an integral control for ISTA by using the hyperparameter as a control input.

The contributions of the paper are twofold. Firstly, we develop the proposed approach, and we analyze its convergence for strongly convex cost functions. Secondly, we illustrate some numerical results to compare the proposed method and state-of-the-art gradient-based techniques in strongly convex and non-strongly convex Lasso problems.

We organize the paper as follows. In Sec. \ref{sec:PS}, we state the problem. In Sec. \ref{sec:PA}, we develop the proposed approach and analyze its convergence in Sec. \ref{sec:CA}. Then, we validate and extend the theoretical results through numerical simulations in Sec. \ref{sec:NR}. Finally, we draw some conclusions.
\section{Problem Statement}\label{sec:PS}
In this work, we consider optimization problems of the kind
$\min_{x\in\R^n}f(x)$
where $f :\R^n \mapsto \R^+ $ is convex, differentiable and admits a sparse minimizer that we aim to estimate.

Since the proposed study potentially addresses high-dimensional data problems, solving $\nabla f(x)=0$ is not a viable way to estimate the sparse minimizer. Moreover, the problem is not well-posed when $f$ has multiple minimizers. We resort to gradient-based methods for these motivations to achieve the desired solution.

To promote sparsity, we modify the problem into
\begin{equation}\label{p1}
\begin{split}
&\min_{x\in\R^n}f(x)+ \sum_{i=1}^n \lambda_i |x_i|%\\% ~~~\text{ where}\\
%&\fun(x):=\frac{1}{2}\|Ax-y\|_2^2+\run(x)\\
\end{split}
\end{equation}
where  and $\sum_{i=1}^n \lambda_i |x_i|$ is a  weighted $\ell_1$ regularization with $\lambda_i \geq 0$ for each $i=1,\dots,n$.

If $f$ admits a unique (sparse) minimizer, in principle, the sparsity-promoting regularization term is unnecessary and we can obtain the solution through gradient descent. However, this approach can be very slow; regularization improves the convergence rate at the price of a bias. On the other hand, if $f$ admits multiple minimizers, the $\ell_1$ regularization plays a crucial role to achieve the desired sparse estimate.

To solve \eqref{p1}, we can apply PGM, which iterates a gradient descent step over $f$ with constant stepsize $\tau>0$ and a proximal mapping with respect to the non-smooth regularization. In the case of $\ell_1$ regularization, the proximal operator corresponds to the shrinkage-thresholding operator $\soft_{\tau\lambda}:\R^n\mapsto \R^n$, which is defined componentwise as follows:
\begin{equation}\label{soft_def}
\begin{split}
\soft_{\tau \lambda_i}(z_i)&:=\argmin{x_i\in\R}\left[\tau \lambda_i |x_i|+\frac{1}{2}(x_i-z_i)^2\right],~~z\in\R\\
&= \left\{\begin{array}{ll}
                  z_i-\sign(z_i)\tau\lambda_i&~\text{ if } |z_i|> \tau\lambda_i\\
                  0&~\text{ otherwise.}
                 \end{array}
\right.
\end{split}
\end{equation}
In conclusion, PGM for \eqref{p1}, also known as ISTA, reads as follows: for $k=0,1,2,\dots$,
%\footnote{In some works, the name ISTA refers to PGM applied to Lasso, which is an instance of the class of problems \eqref{p1}}. ISTA reads as follows: for $k=0,1,2,\dots$,
\begin{equation}\label{ISTA}
\begin{split}
&x(k+1)=\soft_{\tau\lambda}\left(x(k)-\tau \nabla f(x(k))\right).
\end{split}
\end{equation}
When $f$ is strongly convex, the map in \eqref{ISTA} is contractive thanks to the non-expansiveness of $\soft_{\tau\lambda}$, which implies convergence to the minimizer of \eqref{p1}, see, e.g., \cite{for10}. More generally, the convergence of PGM to a minimizer of $\eqref{p1}$ is studied, e.g., in \cite{att10, com11}.
We remark that Lasso, as defined in \cite{tib96}, is an instance of problem \eqref{p1} with $f(x)=\frac{1}{2}\|Ax-y\|_2^2$ where $A\in\R^{m,n}$, $y\in\R^m$, and $\lambda\in\R^n$ has all equal components.

The main goal of this work is to tackle the following problem.
\begin{problem}\label{thep}
Given ISTA as defined in \eqref{ISTA}, we aim at designing a feedback control strategy, acting on $\lambda$ as a control input, that minimizes the solution bias in \eqref{p1}, that is, the distance from the minimum of $f$, while preserving the solution sparsity.
\end{problem}
\subsection{Related literature}
In the literature, two lines of research address the development of control strategies to improve the trajectory of ISTA for Lasso with acceleration purposes. The common idea is to design a time-varying $\lambda$ that improves the dynamics of ISTA.
As noticed in \cite{dau08}, for reasonably small $\lambda>0$, ISTA for Lasso firstly minimizes $f$, then it adjusts the $\ell_1$ norm; see \cite[Fig.1]{dau08}. This trajectory is not optimal because it causes a considerable $\ell_1$ overshoot, whose correction is time consuming. To address this issue, the work \cite{dah12} introduces D-ISTA, that is, an ISTA with geometrically decreasing $\lambda(k)>0$. In other terms, in D-ISTA, $\lambda$  can be interpreted as an open-loop control law. However, the design of a convenient control law is critical, as illustrated in \cite[Theorem 3.2]{dah12}.

In contrast, in \cite{fox23}, the authors consider a feedback control approach for $\lambda$ that originates from the use of a non-convex regularization instead of $\ell_1$ norm. The corresponding PGM, called  AD-ISTA, is usually faster than ISTA and other competitors in Lasso problems. %Contextually, based on \cite{bec09}, a fast version called AD-FISTA is introduced. We remark that the control law in AD-ISTA depends on the formulation of the optimization problem, and it does not derive from feedback control theory arguments.

These two approaches enjoy an excellent interpretation in the framework of control methods, and the corresponding control laws derive from sparse optimization considerations. Moreover, their focus is on accelerating ISTA. In contrast, in Problem \ref{thep}, we start from a control perspective and we focus on the bias regulation.
\section{Proposed approach}\label{sec:PA}
As stated in Problem \ref{thep}, the goal of this work is to remove the bias without affecting the solution sparsity. Since $f$ is differentiable and convex by assumption, removing the bias corresponds to achieving $\nabla f =0$.
The key idea of the proposed approach is to regulate $y(k)=\nabla f(x(k))$ to zero by a suitable feedback control law on ISTA, with $\lambda$ as a control input.
A straightforward choice to regulate the output of a dynamical system to zero is to implement an integral control $\lambda(t)=k_i\int_0^t y(k)$, where $k_i\in\R$ is a design parameter. By recasting the integral law into a discrete-time setting, we obtain $\lambda(k+1)=\lambda(k) + k_i y(k)$, where $k_i$ accounts also for the discretization step.

More precisely, the proposed ISTA with integral control, denoted by I-ISTA, is as follows:

\begin{equation}\label{I-ISTA}
\left\{\begin{split}
&x(k+1)=\soft_{\tau\lambda}\left(x(k)-\tau \nabla f(x(k))\right)\\
&\lambda(k+1)=(1-\alpha)\lambda(k) + k_i\nabla f(x(k))
\end{split}\right.
\end{equation}
where $\alpha\in\ (0,1)$ is a correction term useful for the convergence analysis presented in Sec. \ref{sec:CA}.

%\begin{remark}
We notice that the computational complexity of I-ISTA is similar to that of ISTA because the two algorithms differ only for the update of $\lambda$ in \eqref{I-ISTA}. Regarding the storage requirements, I-ISTA requires to save $2n$ variables instead of $n$ at each iteration.
%\end{remark}

%
\section{Convergence analysis}\label{sec:CA}
In this section, we characterize the equilibrium point of \eqref{I-ISTA} in a strongly convex framework and we prove the convergence of I-ISTA.

We consider the following assumptions.
\begin{assumption}\label{ass:mu}
 $f$ is differentiable and $\mu$-strongly convex, i.e., there exists $\mu>0$ such that $f(x)-\frac{\mu}{2}\|x\|_2^2$ is convex.
\end{assumption}
This implies that, for any $x,z\in\R^n$,
\begin{equation}
f(x)\geq f(z) +\nabla f(z)^\top (x-z)+\frac{\mu}{2}\|x-z\|_2^2.
\end{equation}
\begin{assumption}\label{ass:beta}
$f$ is $\beta$-smooth, that is, there exists $\beta>0$ such that $\frac{\beta}{2}\|x\|_2^2-f(x)$ is convex.
\end{assumption}
$\beta$-smoothness corresponds to the $\beta$-Lipschitz continuity of $\nabla f$. Thus, for any $x,z\in\R^n$
\begin{equation}
f(x)\leq f(z)+\nabla f(z)^\top (x-z)+\frac{\beta}{2}\|x-z\|_2^2.
\end{equation}

We refer the reader to, e.g., \cite{cvx} for details.

%We notice that Assumption \ref{ass:mu} may seem too restrictive:
%As discussed in Sec. \ref{sec:PA} if $f$ is $\mu$-strongly convex, then in particular it has a unique minimizer.

The following result states that the equilibrium point of system \eqref{I-ISTA} provides an unbiased solution.
\begin{lemma}\label{lem:1}
Let Assumption \ref{ass:mu} holds. If $\alpha>|k_i|,$ the equilibrium point $(\xstar,\lstar)$ of \eqref{I-ISTA} satisfies $\nabla f(\xstar)=0$ and $\lstar=0$. In particular, $\xstar$ is the unique, hence sparse minimizer of $f$.
\end{lemma}
\begin{proof}
We compute the equilibrium points of \eqref{I-ISTA}.
\begin{equation}\label{ss}
\left\{\begin{split}
&\xstar=\soft_{\tau\lstar}\left(\xstar-\tau \nabla f(\xstar)\right)\\
&0=-\alpha\lstar + k_i\nabla f(\xstar)
\end{split}\right.
\end{equation}
From the second equation, we have
\begin{equation}\label{grad_condition}
\nabla f(\xstar)=\frac{\alpha}{k_i}\lstar.
\end{equation}
By replacing \eqref{grad_condition} in the first equation of \eqref{ss},
\begin{equation}
\xstar=\soft_{\tau\lstar}\left(\xstar-\tau\frac{\alpha}{k_i}\lstar\right).
\end{equation}
Now, for each $j\in\{1,\dots,n\}$ such that $\xstar_j \neq 0$,
\begin{equation}
\xstar_j=\xstar_j-\tau\frac{\alpha}{k_i}\lstar_j-\sign\left(\xstar_j-\tau\frac{\alpha}{k_i}\lstar_j\right)\tau\lstar_j.
\end{equation}
Therefore, $\lstar_j=0$ if $\frac{\alpha}{|k_i|}\neq 1$.
On the other hand, for each $j\in\{1,\dots,n\}$ such that $\xstar_j = 0$,
\begin{equation}\label{l_cond}
\left|\tau\frac{\alpha}{k_i}\lstar_j\right|\leq \tau\lstar_j.
\end{equation}
If $\frac{\alpha}{|k_i|}>1$, then \eqref{l_cond} holds if $\lstar_j=0$.

In conclusion, $\lstar=0$, which implies $\nabla(f(\xstar))=0$ from \eqref{grad_condition}.
Moreover, $\xstar$ is the unique minimizer of $f$ since $\nabla(f(\xstar))=0$ and $f$ is $\mu$-strongly convex.
\end{proof}

Now, let us analyze the convergence of I-ISTA to the equilibrium point described in Lemma \ref{lem:1}.
\begin{proposition}\label{prop}
Let us set $\tau<\frac{2}{\beta}$ and let
\begin{equation}\label{contraction_coeff}
\xi^2=\max\{\sigma^2+k_i^2\beta^2,\tau^2+(1-\alpha)^2\}<\frac{1}{2}.
\end{equation}
If Assumptions \ref{ass:mu} and \ref{ass:beta} hold, each step of I-ISTA is a contractive map, and I-ISTA converges to $(\xstar,\lstar)$.
\end{proposition}
\begin{proof}
As illustrated in \cite{has21}, since $\soft$ is non-expansive, given assumptions \ref{ass:mu}-\ref{ass:beta} for any $x,z\in\R^n$ and $\lambda \in \R_+^n$
\begin{equation}\label{contr1}
\begin{split}
 &\left\|\soft_{\tau\lambda}(x-\tau\nabla f(x))-\soft_{\tau\lambda}(z-\tau\nabla f(z))\right\|_2^2\leq \sigma^2 \|x-z\|_2^2
\end{split}
\end{equation}
where
$
\sigma^2= \max\{(1-\tau\mu)^2,(1-\tau\beta)^2\}.
$
Since $\mu\leq \beta$, if $\tau<\frac{2}{\beta}$ then $\sigma^2\in (0,1)$, that is, $h(x)=\soft_{\tau\lambda}(x-\tau\nabla f(x))$ is contractive; see \cite[Proof of Lemma 1]{has21}.

On the other hand, for any $\lambda,\gamma\in\R_+^n$ and $x\in\R^n$
\begin{equation}\label{contr2}
\begin{split}
 &\left\|\soft_{\tau\lambda}(x)-\soft_{\tau\gamma}(x)\right\|_2^2\leq \tau^2 \|\lambda-\gamma\|_2^2.
\end{split}
\end{equation}
In fact, through componentwise analysis,
\begin{enumerate}
 \item if $|x_i|>\tau\lambda_i>\tau\gamma_i$, then
$\soft_{\tau\lambda_i}(x_i)-\soft_{\tau\gamma_i}(x_i)=\sign(x_i)\tau(\gamma_i-\lambda_i)$;
\item if $\tau\lambda_i>|x_i|>\tau\gamma_i$, then $|\soft_{\tau\lambda_i}(x_i)-\soft_{\tau\gamma_i}(x_i)|=|x_i-\sign(x_i)\tau\gamma_i| \leq \tau|\lambda_i-\gamma_i|$.
\end{enumerate}

By using \eqref{contr1} and \eqref{contr2} and the fact that $(a+b)^2\leq 2a^2+2b^2$, for any $(x,\lambda)\in\R^{2n}$ and $(z,\gamma)\in\R^{2n}$, we conclude
\begin{equation}\label{C1}
\begin{split}
 &\left\|\soft_{\tau\lambda}(x-\tau\nabla f(x))-\soft_{\tau\gamma}(z-\tau\nabla f(z))\right\|_2^2\leq\\
  &2\left\|\soft_{\tau\lambda}(x-\tau\nabla f(x))-\soft_{\tau\lambda}(z-\tau\nabla f(z))\right\|_2^2+\\
 &2\left\|\soft_{\tau\lambda}(z-\tau\nabla f(z))-\soft_{\tau\gamma}(z-\tau\nabla f(z))\right\|_2^2\leq\\
&2\sigma^2\|x-z\|_2^2+2\tau^2\|\lambda-\gamma\|_2^2.
\end{split}
\end{equation}
As to the first equation of \eqref{I-ISTA}, the bound \eqref{C1} implies that
\begin{equation}\label{CC1}
\left\|x(k+1)-\xstar\right\|_2^2\leq 2\sigma^2\|x(k)-\xstar\|_2^2+2\tau^2\|\lambda(k)-\lstar\|_2^2.
\end{equation}

Furthermore,  regarding the second equation of \eqref{I-ISTA}, we notice that for any $(x,\lambda)\in\R^{2n}$ and $(z,\gamma)\in\R^{2n}$
\begin{equation}\label{C2}
\begin{split}
 &\left\|(1-\alpha)\lambda+k_i\nabla f(x)-(1-\alpha)\gamma-k_i\nabla f(z)\right\|_2^2 \\
 &\leq 2(1-\alpha)^2\|\lambda-\gamma\|_2^2+2k_i^2\|\nabla f(x)-\nabla f(z)\|_2^2\\
 &\leq 2(1-\alpha)^2\|\lambda-\gamma\|_2^2+2k_i^2\beta^2\|x-z\|_2^2
\end{split}
\end{equation}
where in the last step we exploit the $\beta$-smoothness of $f$.
Thus,
\begin{equation}\label{CC2}
\begin{split}
 &\left\|\lambda(k+1)-\lstar\right\|_2^2 \\
 &\leq 2(1-\alpha)^2\|\lambda(k)-\lstar\|_2^2+2k_i^2\beta^2\|x(k)-\xstar\|_2^2.
\end{split}
\end{equation}
By summing \eqref{CC1} and \eqref{CC2}, we obtain
\begin{equation}\label{CC3}
\begin{split}
&\left\|x(k+1)-\xstar\right\|_2^2+\left\|\lambda(k+1)-\lstar\right\|_2^2 \\
&\leq 2\xi^2\left(\|x(k)-\xstar\|_2^2+\|\lambda(k)-\lstar\|_2^2\right)
\end{split}
\end{equation}
where $\xi^2=\max\{\sigma^2+k_i^2\beta^2,\tau^2+(1-\alpha)^2\}<\frac{1}{2}$.
This proves that the mapping between $(x(k),\lambda(k))$ and $(x(k+1),\lambda(k+1))$ is contractive with coefficient $2\xi^2$ , thus I-ISTA converges to $(\xstar, \lstar)$ thanks to the Banach fixed-point theorem \cite{ban22}.
\end{proof}

\begin{remark}
The sufficient conditions of Proposition \ref{prop} are quite restrictive: to guarantee the contractivity of I-ISTA, we exploit the $\mu$-strong convexity of $f$ and the bound \eqref{contraction_coeff}, which limits the values of $\alpha$. However, the obtained bounds are not tight, and numerical results prove that the conditions can be relaxed; see Sec. \ref{sec:NR}.
\end{remark}

\section{Numerical results}\label{sec:NR}
In this section, we present some numerical results that support the theoretical convergence results and extend them to non-strongly convex problems.
Moreover, they provide more insight into the trajectory and convergence speed of I-ISTA with respect to state-of-the-art gradient-based approaches.

Specifically, we compare I-ISTA to ISTA and its fast version, FISTA, introduced by \cite{bec09}, and to AD-ISTA and its fast version, AD-FISTA; see \cite{fox23}. As shown in \cite{fox23}, in some Lasso problems, AD-FISTA is the fastest algorithm among state-of-the-art iterative sparse optimization algorithms.
Moreover, we show the behavior of the gradient descent (without sparsity promoting terms) as a benchmark.

Our experiments focus on recovering sparse vectors from linear measurements via Lasso. We consider $\xtrue\in\R^n$ with $n=200$ and sparsity $\|\xtrue\|_0=10\ll n$. We randomly generate the non-zero components of $\xtrue$ through a uniform distribution, with magnitude in $(1,2)$. We aim to recover $\xtrue$ from $y=A\xtrue$, where $A\in\R^{m,n}$ has components independently generated with Gaussian distribution $\mathcal{N}(0,\frac{1}{m})$.
The cost function is $f(x)=\frac{1}{2}\|Ax-y\|_2^2$. We envisage either strongly convex ($m=210>n$) and non-strongly convex ($m=150<n$) cases.

To implement ISTA and FISTA, we consider Lasso with $\lambda_i=10^{-3}$ for each $i=1,\dots,n$. For AD-ISTA and AD-FISTA, we consider a Log-Lasso with initial $\lambda_i=3\times 10^{-3}$ for each $i=1,\dots,n$ and $\epsi = 10^{-2}$.
The gradient stepsize is $\tau=\|A\|_2^{-2}$ for all the algorithms.
For I-ISTA, we set $k_i=10^{-3}$, while $\alpha=0.05$ for $m=210$ and $\alpha=0.02$ for $m=150$.
For all the algorithms, the stop criterion is $\|x(k+1)-x(k)\|_2<10^{-10}$, with a maximum of iterations set to $5\times 10^4$.
\subsection{Strongly convex case: $m=210>n$}
\graphicspath{{FINAL_SIMS_ECC/graphs_210/}}
\begin{figure*}[h!]
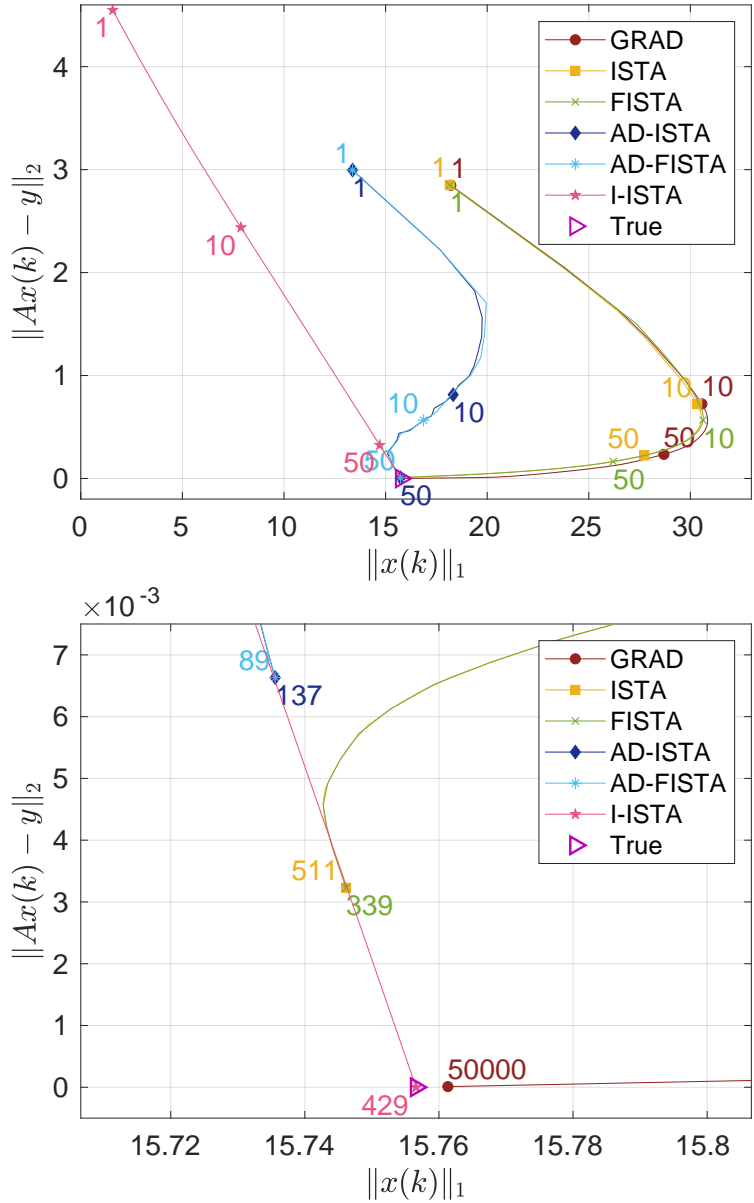

\begin{center}
\includegraphics[width=0.78\columnwidth]{res_vs_l1.pdf}\hspace{1.5cm}
\includegraphics[width=0.78\columnwidth]{res_vs_l1_zoom.pdf}
\caption{Example 1: $m=210$. Residual $\|Ax(k)-y\|_2$ with respect to $\|x(k)\|_1$ in a single run. The curves are parametrized with time. ``True'' refers to the value of $\xtrue$. On the left, we show the overall trajectory; we label iterations 1, 10, 50. On the right,  we magnify the figure around $\xtrue$ and report the convergence step. The gradient descent (GRAD) reaches $\xtrue$, but with a number of iterations larger than the set maximum $5\times10^4$.}
\label{f:1}
\end{center}
\end{figure*}
In Fig. \ref{f:1}, we show the trajectories of all the considered algorithms in the plane $\|Ax(k)-y\|_2$ vs $\|x(k)\|_1$. In contrast to existing ISTA-based approaches, the proposed I-ISTA converges to the true vector $\xtrue$ without bias. The required number of iterations is comparable to ISTA and FISTA. We remark that the number of iterations is a valuable performance metric because the computational complexity of each iteration step is comparable for all the considered algorithms. The gradient descent (GRAD) is also unbiased, although the convergence is very slow.

A remarkable point is the linear trajectory of I-ISTA in the plane $\|Ax(k)-y\|_2$ vs $\|x(k)\|_1$, which highlights an optimal balance between decreasing the residual and increasing the $\ell_1$ norm.
\begin{figure*}[ht!]
\begin{center}
\includegraphics[width=0.78\columnwidth]{rel_error.pdf}\hspace{1.5cm}
\includegraphics[width=0.78\columnwidth]{res.pdf}
\caption{Example 1: $m=210$. Evolution  of the relative error $\|x(k)-\xtrue\|_2/\|\xtrue\|_2$ (left) and of the residual $\|Ax(k)-y\|_2$ (right), averaged over 100 runs.}
\label{f:2}
\end{center}
\end{figure*}

To substantiate these findings, in Fig. \ref{f:1}, we illustrate the time evolution of the relative error $\|x(k)-\xtrue\|_2/\|\xtrue\|_2$ and of the residual $\|Ax(k)-y\|_2$, averaged over 100 different runs.
\begin{figure*}[ht!]
\begin{center}
\includegraphics[width=0.78\columnwidth]{support_error.pdf}\hspace{1.5cm}
\includegraphics[width=0.78\columnwidth]{l0.pdf}
\caption{Example 1: $m=210$. Evolution of the support error $\sum_{i=1}^n |\I(x_i(k)-\I(\xtrue_i)|$ (left) and of the sparsity level $\|x(k)\|_0$ (right), averaged over 100 runs. The graphs on the support error are interrupted when the error is null.}
\label{f:3}
\end{center}
\end{figure*}
Finally, we show the evolution of the estimated support, i.e., the set of non-zero components, in Fig. \ref{f:3}. We notice that I-ISTA identifies the correct support; in fact, the support error defined as $\sum_{i=1}^n |\I(x_i(k)-\I(\xtrue_i)|$, where $\I$ denotes the indicator function $\I(z)=\|z\|_0$ for $z\in \R$, goes to zero in all the runs. Moreover, I-ISTA identifies the correct support after a number of iterations  comparable to AD-ISTA and AD-FISTA. In Fig. \ref{f:3} (right), we highlight a peculiarity of I-ISTA: in contrast to the competitors, it builds the support from below, without the usual ``false positives'' phase that characterizes the ISTA-based methods. This feature can be of interest for all those applications where the transient overestimated support can cause serious false alarms; this is the case, for example, of secure state estimation problems in cyber-physical systems, where the support denotes a subset of sensors under malicious attacks; see, e.g., \cite{fox23sysid} for details.

\subsection{Non-strongly convex case: $m=150<n$}
In this section, we duplicate the numerical simulations with $m=150$. The primary outcome is that I-ISTA converges also for non-strongly convex problems, extending the landscape concerning the proposed convergence analysis.

\begin{figure*}[h!]
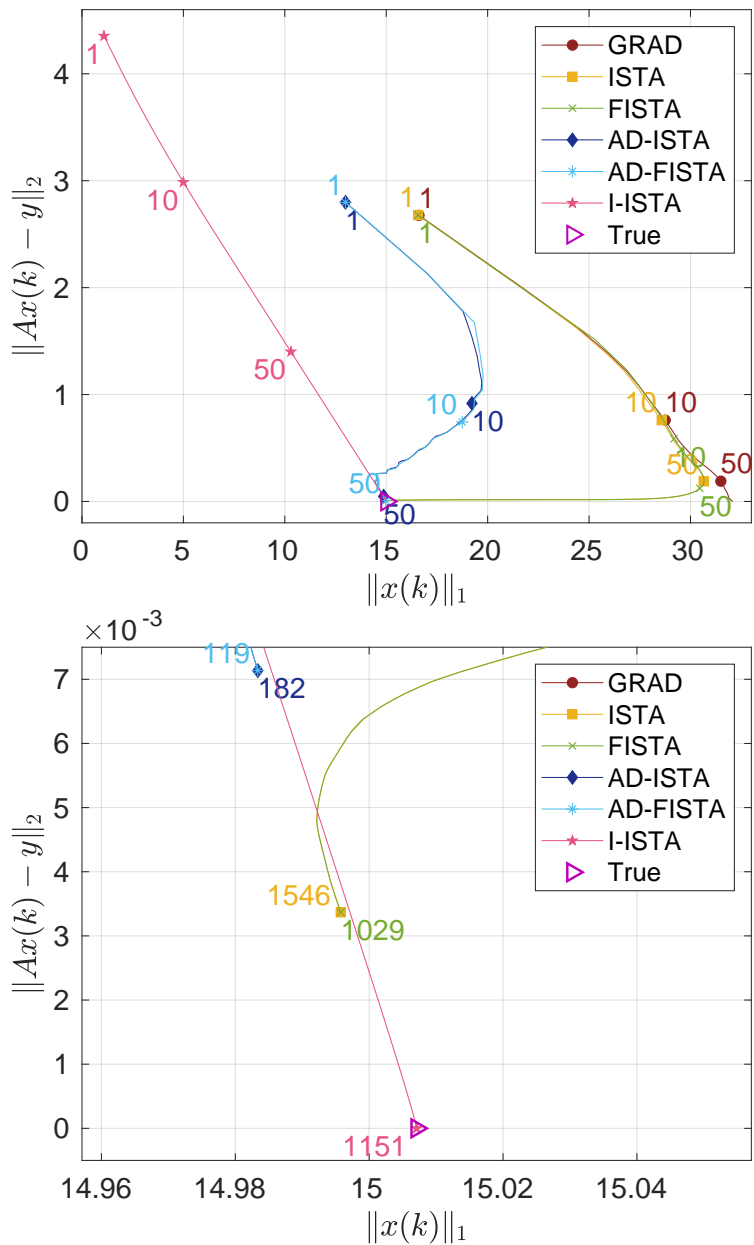

\begin{center}
\includegraphics[width=0.78\columnwidth]{res_vs_l1.pdf}\hspace{1.5cm}
\includegraphics[width=0.78\columnwidth]{res_vs_l1_zoom.pdf}
\caption{Example 2: $m=150$. Residual $\|Ax(k)-y\|_2$ with respect to $\|x(k)\|_1$ in a single run. The curves are parametrized with time. ``True'' refers to the value of $\xtrue$. On the left, we show the overall trajectory; we label iterations 1, 10, 50. On the right,  we magnify the figure around $\xtrue$ and report the convergence step.}
\label{f:11}
\end{center}
\end{figure*}

In this case, $f$ has infinitely many minimizers, and the gradient descent is ineffective because it does not converge to the sparse solution. In contrast, I-ISTA converges precisely to the desired solution, as we can see in Fig. \ref{f:11} and \ref{f:22}. As in the strongly convex case, the convergence time of I-ISTA is comparable to ISTA and FISTA, while the support stabilization time, depicted in Fig. \ref{f:33}, is comparable to AD-ISTA and AD-FISTA.

\graphicspath{{FINAL_SIMS_ECC/graphs_150/}}
%
% In Fig. \ref{f:1}, we show the trajectories of I-ISTA and competitors in the plane $\|Ax(k)-y\|_2$ vs $\|x(k)\|_1$. In contrast to ISTA-based approaches (i.e., ISTA, FISTA, AD-ISTA, and AD-FISTA), the proposed I-ISTA converges to the true vector $\xtrue$, without bias. The number of iterations to converge is comparable to ISTA and FISTA. We remark that the number of iterations is a valuable performance metric because the computational complexity of each iteration step is analogous for all the considered algorithms. The gradient descent (GRAD) is also unbiased, but the convergence is very slow.
%
% A remarkable point is the linear trajectory of I-ISTA in the plane $\|Ax(k)-y\|_2$ vs $\|x(k)\|_1$, which highlights an optimal balance between decreasing the residual and increasing the $\ell_1$ norm.

\begin{figure*}[ht!]
\begin{center}
\includegraphics[width=0.78\columnwidth]{rel_error.pdf}\hspace{1.5cm}
\includegraphics[width=0.78\columnwidth]{res.pdf}
\caption{Example 2: $m=150$. Evolution  of the relative error $\|x(k)-\xtrue\|_2/\|\xtrue\|_2$ (left) and the residual $\|Ax(k)-y\|_2$ (right), averaged over 100 runs.}
\label{f:22}
\end{center}
\end{figure*}
%
% To substantite this finding, in Fig. \ref{f:1}, we illustrate the time evolution of the relative error $\|x(k)-\xtrue\|_2/\|\xtrue\|_2$and of the residual $\|Ax(k)-y\|_2$, averaged over 100 different runs.

\begin{figure*}[ht!]
\begin{center}
\includegraphics[width=0.78\columnwidth]{support_error.pdf}\hspace{1.5cm}
\includegraphics[width=0.78\columnwidth]{l0.pdf}
\caption{Example 2: $m=150$. Evolution of the support error $\sum_{i=1}^n |\I(x_i(k)-\I(\xtrue_i)|$ (left) and of the sparsity level $\|x(k)\|_0$ (right), averaged over 100 runs. The graphs on the support error are interrupted when the error is null.}
\label{f:33}
\end{center}
\end{figure*}
% Finally we show the evolution of the estimated support, i.e., the set of non-zero components, in Fig. \ref{f:3}. Firstly, we notice that I-ISTA identifies the correct support; in fact, the support error defined as $\sum_{i=1}^n |\I(x_i(k)-\I(\xtrue_i)|$, where $\I$ denotes the indicator funcion $\I(z)=\|z\|_0$ for $z\in \R$, goes to zero in all the 100 runs. Moreover, we notice that I-ISTA identifies the correct support with a convergence speed comparable to AD-ISTA and AD-FISTA. In Fig. \ref{f:3} (right), we highlight a peculiarity of I-ISTA: in contrast to the competitors, it builds the support from below, without the usual phase of ``false positives'' that characterizes the ISTA-based methods. This feature can be of interest for all those applications where the transient overestimated support can cause serious false alarms; this is the case, for example, of secure state estimation problems in cyber-physical systems, where the support denotes a subset of sensors under malicious attacks; see, e.g., \cite{fox23sysid} for details.
%
\begin{table}[h]
\renewcommand{\arraystretch}{1.0}
\centering
{\footnotesize{
    \begin{tabular}{r| c c | c c}
     & $m=210$ & & $m=150$ & \\
    \hline
    Algorithm & Conv. & Supp. stab. & Conv. & Supp. stab.\\
    \hline
  GRAD     &     47183.57&           8002.83 &       2687.42 &     --\\
  ISTA     &            486.36&             382.36&   1761.47&      1617.16\\
  FISTA    &            322.40&             255.76&   1172.71&      1079.11\\
  AD-ISTA  &            123.80&              23.70&    172.80&        45.07\\
  AD-FISTA &             80.01&              16.70&    113.17&        31.09\\
  I-ISTA   &            426.33&               8.23&   1107.80&        25.40\\
  \end{tabular}
    }}
    \caption{Mean number of iterations to converge and to stabilize the support, over 100 random runs.}
    \label{tab:1}
\end{table}
In Table \ref{tab:1}, we summarize the data about convergence and support stabilization times.
\section{Conclusions}\label{sec:C}
In this work, we propose and analyze I-ISTA, a proximal gradient method for $\ell_1$-regularized sparse optimization problems with an integral control that removes the bias. We analyze the convergence of the algorithm in the framework of strongly convex problems, while numerical results extend its validity to non-strongly convex problems. The ``linear'' trajectory of I-ISTA yields a fast stabilization of the support estimate without support overshoot. Significant extensions under investigation are the convergence analysis in non-strongly convex frameworks and the robustness to noise. Furthermore, we are studying how to expand and refine the approach by considering controllers that are more sophisticated than the integral one.

% AD-ISTA, a variant of ISTA developed by applying the proximal gradient method to Log-Lasso. AD-ISTA converges in less iterations with respect to ISTA, FISTA and ADMM, thanks to an adaptive shrinkage hyperparameter, that limits the increase of the $\ell_1$-norm during the first phase. Moreover, by applying the principles of FISTA, we also propose the accelerated version AD-FISTA. Through numerical experiments, we verify that AD-ISTA is faster than the state-of-the-art algorithms for Lasso and that we obtain a further acceleration with AD-FISTA. Possible extensions of this work include the rigorous proof of the convergence rate and the generalization to sparse optimization problems different from Lasso.
% %

\bibliographystyle{IEEEtran}
\bibliography{refs}             % bib file to produce the bibliography
\end{document}